\documentclass[onefignum,onetabnum]{siamart171218}
\usepackage{bbm}
\usepackage{subfig}
\usepackage{amsfonts}
\usepackage{graphicx}
\usepackage{epstopdf}
\usepackage{algorithm}
\usepackage[noend]{algpseudocode}
\usepackage{ulem}
\DeclareMathOperator*{\argmin}{arg\,min}

\DeclareMathOperator{\arcsinh}{arcsinh}

\ifpdf
  \DeclareGraphicsExtensions{.eps,.pdf,.png,.jpg}
\else
  \DeclareGraphicsExtensions{.eps}
\fi

\newsiamremark{remark}{Remark}
\newsiamremark{hypothesis}{Hypothesis}
\crefname{hypothesis}{Hypothesis}{Hypotheses}
\newsiamthm{claim}{Claim}
\crefname{claim}{Claim}{Claim}

\headers{TD}{TD}

\title{Second Order Threshold Dynamics Schemes for\\ Two Phase Motion by Mean Curvature}

\author{Alexander Zaitzeff\thanks{Department of Mathematics, University of Michigan, Ann Arbor, MI 48109, USA
  (\email{azaitzef@umich.edu}, \email{esedoglu@umich.edu}).}
\and Selim Esedo\=glu\footnotemark[1]
\and Krishna Garikipati\thanks{Departments of Mechanical Engineering, and Mathematics, Michigan Institute for Computational Discovery \& Engineering, University of Michigan, Ann Arbor, MI 48109, USA (\email{krishna@umich.edu}).}}

\usepackage{amsopn}

\DeclareMathOperator\erf{erf}

\ifpdf
\hypersetup{
  pdftitle={Second Order TD},
  pdfauthor={ZEG}
}
\fi

\begin{document}

\maketitle

\begin{abstract}
The threshold dynamics algorithm of Merriman, Bence, and Osher is only first order accurate in the two-phase setting.
Its accuracy degrades further to half order in the multi-phase setting, a shortcoming it has in common with other related, more recent algorithms such as the equal surface tension version of the Voronoi implicit interface method.
As a first, rigorous step in addressing this shortcoming, we present two different second order accurate versions of two-phase threshold dynamics.
Unlike in previous efforts in this direction, we present careful consistency calculations for both of our algorithms.
The first algorithm is consistent with its limit (motion by mean curvature) up to second order in any space dimension.
The second achieves second order accuracy only in dimension two, but comes with a rigorous stability guarantee (unconditional energy stability) in any dimension -- a first for high order schemes of its type.
\end{abstract}



\section{Introduction}
In this paper, we will describe new, second order accurate in time versions of a popular algorithm for simulating the motion of interfaces by mean curvature known as threshold dynamics.
The original version of the algorithm, which is only first order accurate in time in the two-phase setting, was proposed by Merriman, Bence, and Osher in \cite{mbo92, mbo94}.
Since then, many extensions of the algorithm have been given, for instance to multiphase mean curvature motion, where it has proven particularly useful and flexible.
There have also been high order accurate versions of the algorithm proposed in several previous studies, discussed in detail in Section \ref{sec:pre}

For a $(d-1)$-dimensional smooth interface $\Gamma\subset\mathbb{R}^d$ given as the boundary of a set $\Sigma\subset\mathbb{R}^d$, the original threshold dynamics algorithm generates a discrete in time approximation to its motion by mean curvature as follows:
\begin{algorithm}[H]
\caption{Original Threshold Dynamics of MBO'92}
\label{alg:mbo}
Fix a time step size $\delta t>0$. Alternate the following steps:
\begin{enumerate}
\item Convolution:
\begin{equation*}
\phi(x) = G_{\delta t} * \mathbf{1}_{\Sigma^k}.
\end{equation*}
\item Thresholding:
\begin{equation*}
\Sigma^{k+1} = \Big\{ x \, : \, \phi(x) \geq \frac{1}{2} \Big\}.
\end{equation*}
\end{enumerate}
\end{algorithm}	
\noindent Here, $G_{\delta t}$ is the Gaussian kernel:
\begin{equation*}
G_{ \delta t}(\mathbf{x})=\frac{1}{(4 \pi \delta t)^{d/2}}e^{-\frac{|\mathbf{x}|}{4\delta t}}
\end{equation*}

Our goal in this paper is to take a step towards providing more accurate versions of threshold dynamics.
The accuracy issue is particularly acute in the multi-phase setting, where it decreases to half-order in time due to the presence of junctions.
Here, we focus on the easier yet still challenging two-phase setting, to find a version of Algorithm \ref{alg:mbo} that
\begin{itemize}
\item Maintains the simplicity and spirit of the original threshold dynamics Algorithm \ref{alg:mbo},
\item Achieves second order accuracy in time,
\item Maintains the variational interpretation, and the resulting stability properties, given in \cite{eo} for the original Algorithm \ref{alg:mbo}.
\end{itemize}

The paper is organized as follows:
\begin{itemize}
\item In \cref{sec:pre}, we recall previous efforts in designing second order versions of threshold dynamics.
\item In \cref{sec:ker}, we discuss necessary conditions for second order accuracy. 
\item In \cref{sec:naturalextrap}, we present our first new algorithm: a natural two kernel extrapolation method, applied to the original threshold dynamics algorithm, to achieve second order accuracy in any space dimension.
\item In \cref{sec:ms}, we present our second new algorithm: a multi-step method that is second order accurate in two space dimensions, and unconditionally energy stable in any dimension.
\item In \cref{sec:numerical}, we provide numerical verification of the advertised order of accuracy for both of our new algorithms.
\end{itemize}

The code for \cref{sec:numerical} is publicly available, and can be found at \url{https://github.com/AZaitzeff/secondorderTD}.

\section{Previous Work}
\label{sec:pre}
In \cite{ruuth1998efficient}, Ruuth proposed the following method based on Richardson extrapolation to jack up the order of accuracy of Algorithm \ref{alg:mbo} to second order in time:
\begin{algorithm}[H]
\caption{Ruuth's Second Order Threshold Dynamics}
\label{alg:ruuth}
Fix a time step size $\delta t>0$.
Set $\phi^k(x) = \mathbf{1}_{\Sigma^0}(x)$.
Alternate the following steps:
\begin{enumerate}
\item First half time step:
$$ \Sigma_1 = \Big\{ x \, : \, G_{\delta t/2} * \phi^k  \geq \frac{1}{2} \Big\}.$$
\item Second half time step:
$$ \Sigma_2 = \Big\{ x \, : \, G_{\delta t/2} * \mathbf{1}_{\Sigma_1} \geq \frac{1}{2} \Big\}.$$
\item Full time step:
$$ \Sigma_3 = \Big\{ x \, : \, G_{\delta t} * \phi^k  \geq \frac{1}{2} \Big\}.$$
\item Linear combination:
$$ \phi^{k+1} = 2\mathbf{1}_{\Sigma_2} - \mathbf{1}_{\Sigma_3}. $$
\end{enumerate}
\end{algorithm}

Although numerical experiments indicate this version indeed improves the accuracy in time to second order for smooth interfaces undergoing two-phase motion by mean curvature, the algorithm sacrifices an attractive simplicity of the original MBO scheme: it no longer generates binary functions exclusively that can be naturally identified with sets.
Perhaps more importantly, there appears to be no clear extension of the variational interpretation given in \cite{eo} for the original Algorithm \ref{alg:mbo} to this case. No comparison principle is expected to hold, as a non-positive weighted sum is involved.
Hence, there is no rigorous result indicating the stability of the algorithm (or its convergence).

In \cite{gh}, Grzhibovskis \& Heinz propose another approach to improving the order of accuracy of Algorithm \ref{alg:mbo} to second order.
The idea is natural: To cancel out the leading order error in threshold dynamics by taking a linear combination of convolutions with two different radially symmetric kernels:
\begin{algorithm}[H]
\caption{Algorithm of Grzhibovskis \& Heinz}
\label{alg:gh}
\begin{enumerate}
\item Convolution step:
$$ \phi(x) = \Big( \alpha K_1 - \beta K_2 \Big) * \mathbf{1}_{\Sigma^k}.$$
\item Threhsolding step:
$$ \Sigma^{k+1} = \Big\{ x \, : \, \phi(x) \geq \frac{1}{2} \Big\}.$$
\end{enumerate}
\end{algorithm}
\noindent The coefficients $\alpha$ and $\beta$ are chosen so that the leading order correction to curvature in the standard consistency calculation for the original threshold dynamics Algorithm \ref{alg:mbo} cancels out.
Crucially, this necessitates that the resulting combined convolution kernel changes sign, even when the individual kernels $K_1$ and $K_2$ are positive.
This means that the resulting algorithm can violate the comparison principle.
But far more importantly, we show in \cref{sec:ker}, that this algorithm {\it does not give second order accuracy} in time; it merely achieves a more accurate evaluation of the mean curvature term at every time step.
In general, the dynamics generated is still only first order accurate, at least without being much more specific and deliberate about the choice of the kernels $K_1$ and $K_2$ -- which the authors do not specify.
(For example, in case both $K_1$ and $K_2$ are Gaussians -- with potentially different mass and/or width -- no choice of the coefficients $\alpha$ and $\beta$ results in a second or higher order accurate in time scheme for motion by mean curvature.)

In this paper, we will provide truly second order accurate in time versions of Algorithm \ref{alg:mbo} that maintain its elegant and simple nature.
Moreover, we will be able to provide rigorous stability results for our new algorithms.

\section{Second Order Motion by Mean Curvature}
\label{sec:ker}
First, we need to identify how far a surface travels under motion by mean curvature. In the vicinity of a point of interest on the surface, which we take to be the origin, let the surface be given as the graph of a smooth function $f(x,y,t): \mathbb{R}^{2}\times[0,\infty) \to \mathbb{R}$ with $f(0,0,0)=0$, $f_{x}(0,0,0)=0$ and $f_{y}(0,0,0)=0$. Since the normal direction changes during the evolution, it is easier to insist that the numerically generated solution intersects a fixed line at nearly the same location as the true solution, at any given time. Thus, we will calculate how far the surface travels along the $z$-axis under mean curvature motion and under our algorithms. For a surface given as the graph of a function, motion by mean curvature takes the following form:

\begin{equation}
\label{eq:mcmpde}
f_t=\frac{f_{xx}(1+ f_{y}^2)-2f_{x}f_{y}f_{xy}+f_{yy}(1+ f_{x}^2)}{1+f_x^2+f_y^2}.
\end{equation}
By a straightforward Taylor expansion we have for small $t$

    \begin{multline}
    \label{eq:mctaylor}
    f(0,0,t) = t \bigg[ f_{xx}+f_{yy} \bigg]\\+t^2\bigg[\frac{1}{2}(f_{xxxx}+2f_{xxyy}+f_{yyyy})-( f_{xx}^3+3f_{xx}f_{xy}^2+3f_{yy}f_{xy}^2+f_{yy}^3)\bigg]+O(t^3)
    \end{multline}
\\
where the functions on the right hand side are evaluated at $(0,0,0)$. Over the course of this paper, we will denote $f(0,0,0)$ as $f$, $f_x(0,0,0)$ as $f_x$, etc. for convenience.\\

It has been known and verified by Taylor expansion in previous publications (e.g. Ruuth~\cite{ruuth1998efficient}) that standard threshold dynamics is 1st order accurate. We will include the expansion of standard threshold dynamics (\cref{alg:mbo}) here as a simple example of the method we use throughout this paper.  Let $\Sigma^0=\{(x,y,z): z \leq f(x,y,0)\}$.
We work out the convolution of a Gaussian kernel with a characteristic function in \cref{sec:taylor} of the appendix where we had to keep many more terms than in previous works to achieve our goals in this paper. Applying calculation \ref{eq:fconv} to our problem we have:
\begin{align}
\label{eq:onestep}
\begin{split}
&G_t*\mathbf{1}_{\Sigma^0}(0,0,z)-\frac{1}{2}
=-\frac{z}{2\sqrt{\pi t}}+\frac{z^3}{24\sqrt{\pi }t^{3/2}}+\frac{\sqrt{t}}{2\sqrt{\pi}}(f_{xx}+f_{yy})\\&+\frac{t^{3/2}}{4\sqrt{\pi}}(f_{xxxx}+2f_{xxyy}+f_{yyyy})
    -\frac{z^2}{8\sqrt{\pi t}}(f_{xx}+f_{yy})\\&+\frac{z\sqrt{t}}{2\sqrt{\pi}} (\frac{3}{4}f_{xx}^2+\frac{3}{4}f_{yy}^2+\frac{1}{2}f_{xx}f_{yy}+f_{xy}^2)\\&-\frac{t^{3/2}}{2\pi}(\frac{5}{4} f_{xx}^3+\frac{5}{4} f_{yy}^3+\frac{3}{4}f_{xx}f_{yy}^2+\frac{3}{4}f_{xx}^2f_{yy}+3f_{xx}f_{xy}^2+3f_{yy}f_{xy}^2)+\text{h.o.t.}
    \end{split}
\end{align}
Next set \cref{eq:onestep} to zero and solve for $z$ by using the ansatz $z=z_1t+z_2t^2+\text{remainder}$ and matching terms of the same order in $t$. Up to second order:
\begin{multline}
\label{eq:tdtaylor3d}
z=t\bigg[ f_{xx}+f_{yy} \bigg]\\+t^2\bigg[\frac{1}{2}(f_{xxxx}+2f_{xxyy}+f_{yyyy})-\frac{2}{3}( f_{xx}^3+3f_{xx}f_{xy}^2+3f_{yy}f_{xy}^2+f_{yy}^3)\bigg]+O(t^3).
\end{multline}
Equation \cref{eq:tdtaylor3d} gives the location of $\partial \Sigma^{1}$ along the $z$-axis. The equation \cref{eq:tdtaylor3d} matches the two dimensional version calculated by Ruuth~\cite{Ruuth_1996}. Additionally, for three dimensions some of the terms in \cref{eq:onestep} are calculated by Grzhibovskis \& Heinz~\cite{grzhibovskis2008convolution}. Their paper, focusing on Willmore flow, did not require all the terms calculated in \cref{eq:onestep}.\\
Comparing \cref{eq:tdtaylor3d} to the location of the interface under mean curvature motion \cref{eq:mctaylor}, we see that threshold dynamics is only first order accurate in time. At this point, we can also already see why \cref{alg:gh} of Grzhibovskis \& Heintz cannot be second order accurate: It would merely move the surface by $t(f_{xx}+f_{yy}) + O(t^3)$, which would still make it only a first order
accurate approximation of the right hand side of \cref{eq:mctaylor}.
In the next two sections, we will present two second order methods. For each method, we will show that they do indeed match motion by mean curvature \cref{eq:mctaylor} up to second order in the normal direction.

\section{A More Natural Two Kernel Extrapolation}
\label{sec:naturalextrap}
Our first method is a two stage algorithm using two different Gaussian kernels with differing amplitudes and widths. We detail the method in \cref{alg:twostep}.
\begin{algorithm}[H]
\caption{Natural Two Kernel Extrapolation}
\label{alg:twostep}
Fix a time step size $\delta t>0$.
Alternate the following steps:
\begin{enumerate}
\item First stage:
$$ \bar{\Sigma} = \Big\{ x \, : \, G_{\delta t/2} * \mathbf{1}_{\Sigma^k}  \geq \frac{1}{2} \Big\}.$$
\item Second stage:
$$ \Sigma^{k+1} = \Big\{ x \, : \,\frac{\sqrt{2} G_{\delta t/2} * \mathbf{1}_{\bar{\Sigma}}-G_{\delta t} *  \mathbf{1}_{\Sigma^k}}{\sqrt{2}-1} \geq \frac{1}{2} \Big\}.$$
\end{enumerate}
\end{algorithm}
Whether \cref{alg:twostep} is unconditionally stable is currently unknown.  We will devote the rest of the section to showing that \cref{alg:twostep} is indeed second order. 

\subsection{Consistency of \cref{alg:twostep}}
Once again let $\Sigma^0=\{(x,y,z): z \leq f(x,y,0)\}$ for $f$ defined in \cref{sec:ker}.
First, we need to find the location and curvature of $\bar{\Sigma}$ along the $z$-axis.  
Let $h(x,y)$ be defined by the requirement that $G_{t/2}*\mathbf{1}_{\Sigma^0}(x,y,h(x,y))=\frac{1}{2}$, so that $\bar{\Sigma}=\{(x,y,z): z \leq h(x,y)\}$. From \cref{eq:tdtaylor3d} we have that 
\begin{multline}
\label{eq:vbar}
h(0,0)=\frac{t}{2}\bigg[ f_{xx}+f_{yy} \bigg]\\+\frac{t^2}{4}\bigg[\frac{1}{2}(f_{xxxx}+2f_{xxyy}+f_{yyyy})-\frac{2}{3}( f_{xx}^3+3f_{xx}f_{xy}^2+3f_{yy}f_{xy}^2+f_{yy}^3)\bigg]+O(t^3)
\end{multline}
Now we would like to find $h_{xx}(0,0)$ and $h_{yy}(0,0)$.
From \cref{eq:der2} in the appendix, replacing $h$ with $f$, $z$ with $h$ and $t$ with $t/2$, we have that $h_{xx}(0,0)$ satisfies
\small
\begin{align}
\label{eq:vbarderstep}
\begin{split}
0=&- \frac{\sqrt{2}h_{xx}}{2\sqrt{\pi t}}+\frac{\sqrt{2}h^2h_{xx}}{4\sqrt{\pi}t^{3/2}}+\frac{\sqrt{2}f_{xx}}{2\sqrt{\pi t}}+\frac{\sqrt{t}}{2\sqrt{2 \pi}}(f_{xxxx}+f_{xxyy})\\
    &-\frac{\sqrt{2}h^2}{4\sqrt{\pi}t^{3/2}}f_{xx}+\frac{\sqrt{2}h}{2\sqrt{\pi t}} (\frac{3}{2}f_{xx}^2+\frac{1}{2}f_{xx}f_{yy}+f_{xy}^2) \\&-\frac{\sqrt{t}}{2\sqrt{2 \pi}}\bigg(\frac{15}{4} f_{xx}^3+\frac{3}{4}f_{xx}f_{yy}^2+\frac{3}{2}f_{xx}^2f_{yy}+6f_{xx}f_{xy}^2+3f_{yy}f_{xy}^2\bigg)\\
    &+\frac{\sqrt{t}h_{xx}}{2\sqrt{2 \pi}}\bigg(\frac{3}{4}f_{xx}^2+\frac{3}{4}f_{yy}^2+\frac{1}{2}f_{xx}f_{yy}+f_{xy}^2\bigg) -\frac{\sqrt{2}hh_{xx}}{4\sqrt{\pi t}}(f_{xx}+f_{yy}) +\text{h.o.t.}
\end{split}
\end{align}
\normalsize
Plugging in $h$ from \cref{eq:vbar} and solving for $h_{xx}$ in \cref{eq:vbarderstep}, using the ansatz $h_{xx}=h_0+th_1+\text{remainder}$, we have 
\[h_{xx}=f_{xx}+\frac{t}{2}(f_{xxxx}+f_{xxyy}-2f_{xx}^3-4f_{xx}f_{xy}^2-2f_{yy}f_{xy}^2)+O(t^2)
\]
We can use similar steps to find $h_{yy}$. Putting $h_{xx}$ and $h_{yy}$ together we arrive at  
\begin{multline}
\label{eq:vbarder}
    h_{xx}+h_{yy}=f_{xx}+f_{yy}\\+\frac{t}{2}\bigg[(f_{xxxx}+2f_{xxyy}+f_{yyyy})-2(f_{xx}^3+3f_{xx}f_{xy}^2+3f_{yy}f_{xy}^2+f_{yy}^3)\bigg]+O(t^2).
\end{multline}

Now we can solve for the location of $\Sigma^{1}$ along the $z$-axis. From the expansion of the Gaussian kernel convoluted with a characteristic function, \cref{eq:fconv},
we have
\small
\begin{align}
\label{eq:penultstep}
\begin{split}
    [\sqrt{2}&G_{\delta t/2} * \mathbf{1}_{\bar{\Sigma}}-G_{\delta t} *  \mathbf{1}_{\Sigma^0}](0,0,z)-\frac{\sqrt{2}-1}{2}\\
    =&-\frac{z}{\sqrt{\pi t}}+\frac{z^3}{6\sqrt{\pi }t^{3/2}}+\frac{h}{\sqrt{\pi t}}+\frac{\sqrt{t}}{2\sqrt{\pi}}(h_{xx}+h_{yy})\\&+\frac{t^{3/2}}{8\sqrt{\pi}}(h_{xxxx}+2h_{xxyy}+h_{yyyy})
    -\frac{z^2}{2\sqrt{\pi}t^{3/2}}h-\frac{z^2}{4\sqrt{\pi t}}(h_{xx}+h_{yy})+\frac{z}{2\sqrt{\pi }t^{3/2}}h^2\\&+\frac{z}{2 \sqrt{\pi t}}h(h_{xx}+h_{yy})+\frac{z\sqrt{t}}{2\sqrt{\pi}} (\frac{3}{4}h_{xx}^2+\frac{3}{4}h_{yy}^2+\frac{1}{2}h_{xx}h_{yy}+h_{xy}^2)-\frac{h^3}{6\sqrt{\pi}t^{3/2}} \\&-\frac{h^2}{4\sqrt{\pi t}}(h_{xx}+h_{yy})-\frac{\sqrt{t}}{2\sqrt{\pi}}h(\frac{3}{4}h_{xx}^2+\frac{3}{4}h_{yy}^2+\frac{1}{2}h_{xx}h_{yy}+h_{xy}^2)\\&-\frac{t^{3/2}}{4\sqrt{\pi}}(\frac{5}{4} h_{xx}^3+\frac{5}{4} h_{yy}^3+\frac{3}{4}h_{xx}h_{yy}^2+\frac{3}{4}h_{xx}^2h_{yy}+3h_{xx}h_{xy}^2+3h_{yy}h_{xy}^2)\\
    &+\frac{z}{2\sqrt{\pi t}}-\frac{z^3}{24\sqrt{\pi }t^{3/2}}-\frac{\sqrt{t}}{2\sqrt{\pi}}(f_{xx}+f_{yy})\\&-\frac{t^{3/2}}{4\sqrt{\pi}}(f_{xxxx}+2f_{xxyy}+f_{yyyy})
    +\frac{z^2}{8\sqrt{\pi t}}(f_{xx}+f_{yy})\\&-\frac{z\sqrt{t}}{2\sqrt{\pi}} (\frac{3}{4}f_{xx}^2+\frac{3}{4}f_{yy}^2+\frac{1}{2}f_{xx}f_{yy}+f_{xy}^2) \\&+\frac{t^{3/2}}{2\sqrt{\pi}}(\frac{5}{4} f_{xx}^3+\frac{5}{4} f_{yy}^3+\frac{3}{4}f_{xx}f_{yy}^2+\frac{3}{4}f_{xx}^2f_{yy}+3f_{xx}f_{xy}^2+3f_{yy}f_{xy}^2)+\text{h.o.t}
    \end{split}
\end{align}
\normalsize
Note that the derivatives of $h$ match the corresponding derivatives of $f$ up to order $t$ (stated precisely in \cref{eq:blakentder}). Substituting \cref{eq:vbar} for $h$, \cref{eq:vbarder} for $h_{xx}+h_{yy}$ and \cref{eq:blakentder} for the other other derivatives of $h$ in \cref{eq:penultstep}, and simplifying, we have:
\small
\begin{align}
\label{eq:finalstep}
\begin{split}
[\sqrt{2}&G_{\delta t/2} * \mathbf{1}_{\bar{\Sigma}}-G_{\delta t} *  \mathbf{1}_{\Sigma^k}](0,0,z)-\frac{\sqrt{2}-1}{2}\\
    =&-\frac{z}{2\sqrt{\pi t}}+\frac{z^3}{8\sqrt{\pi }t^{3/2}}+\frac{t^{3/2}}{4\sqrt{\pi }}(f_{xxxx}+2f_{xxyy}+f_{yyyy})\\&-\frac{2t^{3/2}}{3\sqrt{\pi }}( f_{xx}^3+3f_{xx}f_{xy}^2+3f_{yy}f_{xy}^2+f_{yy}^3)+\frac{\sqrt{t}}{2\sqrt{\pi}}(f_{xx}+f_{yy})\\&
    -\frac{3z^2}{8\sqrt{\pi t}}(f_{xx}+f_{yy})+\frac{3z\sqrt{t}}{8\sqrt{\pi }}(f_{xx}+f_{yy})^2\\&-\frac{\sqrt{t}}{12\sqrt{\pi}}(f_{xx}+f_{yy})^3-\frac{\sqrt{t}}{4\sqrt{\pi}}(f_{xx}+f_{yy})(\frac{3}{4}f_{xx}^2+\frac{3}{4}f_{yy}^2+\frac{1}{2}f_{xx}f_{yy}+f_{xy}^2)\\&+\frac{t^{3/2}}{4\sqrt{\pi}}(\frac{5}{4} f_{xx}^3+\frac{5}{4} f_{yy}^3+\frac{3}{4}f_{xx}f_{yy}^2+\frac{3}{4}f_{xx}^2f_{yy}+3f_{xx}f_{xy}^2+3f_{yy}f_{xy}^2)+\text{h.o.t.}
    \end{split}
\end{align}
\normalsize
Finally, we set \cref{eq:finalstep} to zero and solve for $z$ up to order $t^2$ to obtain:
\begin{multline}
\label{eq:ms2taylor}
z=t\bigg[ f_{xx}+f_{yy} \bigg]\\+t^2\bigg[\frac{1}{2}(f_{xxxx}+2f_{xxyy}+f_{yyyy})-( f_{xx}^3+3f_{xx}f_{xy}^2+3f_{yy}f_{xy}^2+f_{yy}^3)\bigg]+O(t^3).
\end{multline}
Thus,  up to second order, \cref{alg:twostep} moves the interface in the normal direction by the same amount as mean curvature motion \cref{eq:mctaylor}.\\

The drawback to \cref{alg:twostep} is that we do not know whether it is unconditionally stable. In the next section, we present an algorithm that is unconditionally stable but is only guaranteed to be second order in two dimensions. 

\section{Unconditionally Stable Multistage Methods}
\label{sec:ms}
In this section, we provide a class of unconditionally stable threshold dynamics algorithms that are second order in two dimensions. The original threshold dynamics algorithm (\cref{alg:mbo}) is unconditionally energy stable, specifically:
\begin{equation}
\label{eq:uncond}
    E_t(\Sigma^{k+1})\leq E_t(\Sigma^{k})
\end{equation}
for energy
\begin{equation}
\label{eq:energy}
    E_t(\Sigma)=\int_{\mathbb{R}^2} (1-\mathbf{1}_{\Sigma})G_t*\mathbf{1}_{\Sigma} dxdy.
\end{equation}

Our class of methods preserves property \cref{eq:uncond} while at the same time achieving second order accuracy. We describe our $M$-stage method in \cref{alg:mstage}. 

\begin{algorithm}[H]
\caption{M-Stage Unconditionally Stable Threshold Dynamics}
\label{alg:mstage}
Fix a time step size $\delta t>0$ and a choice of $\gamma$'s such that $\sum_{i=0}^{m-1} \gamma_{m,i} =1$ for $m=1,2,\ldots, M$. Set $\tau=\delta t/\beta_{1,M}$ for $\beta_{1,M}$ defined in \cref{eq:rec} and set $\bar{\Sigma}_0=\Sigma^k$.\\

For $m=1,2,\ldots, M$:

\begin{equation}
    \label{eq:tdms}
    \bar{\Sigma}_m=\Big\{ x:G_\tau*\sum_{i=0}^{m-1}\gamma_{m,i}\mathbf{1}_{\bar{\Sigma}_i} \geq \frac{1}{2} \Big\}.
\end{equation}
Then set $\Sigma^{k+1}=\bar{\Sigma}_M$
\end{algorithm}
Unlike the previous algorithm, \cref{alg:mstage} uses the same kernel at each stage. As will be shown later, this will allow us to prove unconditional stability \cref{eq:uncond}. In the rest of the section, we will derive the consistency equations for $\gamma$, give the conditions on $\gamma$ for unconditional stability to hold, and then give a particular choice of $\gamma$ that makes \cref{alg:mstage} unconditionally stable and second order. We conclude with a discussion of one way to extend \cref{alg:mstage} to higher dimensions. 

\subsection{Consistency Equations}

Similar to in three dimensions, let $f(x,t): \mathbb{R} \times[0,\infty) \to \mathbb{R}$ be a graph that is the interface of a set in $\mathbb{R}^2$ with $f(0,0)=0$, $f_{x}(0,0)=0$.
The distance the graph moves under mean curvature motion along the y-axis is:
\begin{equation}
\label{eq:true}
f(0,t)=tf_{xx}+t^2\bigg[\frac{1}{2}f_{xxxx}-f_{xx}^3\bigg]+O(t^2).
\end{equation}

We now present the consistency equations for \cref{alg:mstage}:

\begin{claim}
\label{claim:cons}
Let $\bar{\Sigma}_i$ be given in \cref{eq:tdms} and let
$$\bar{\Sigma}_0=\Sigma^{k}=\Big\{ x:x\leq f(x,0) \Big\}$$.

Define $h_i$ as
$\bar{\Sigma}_i=\Big\{ x:x\leq h_i(x) \Big\}$, so $h_0(x)=f(x,0)$.
The Taylor expansion of $h_i$ at each stage has the same form as \cref{eq:true}, namely: 
\small
\begin{equation}
\label{eq:hfunc}
h_i(0)=t \beta_{1,i} f_{xx}+t^2(\beta_{2,i}f_{xxxx}-\beta_{3,i}f_{xx}^3)+O(t^3).\end{equation}
    Additionally, the Taylor expansion of the second derivative of $h_i$ has form
    \begin{equation}
\label{eq:hder}
 \frac{d^2}{dx^2} h_i(0)=f_{xx}+t(\beta_{4,i}f_{xxxx}-\beta_{5,i}f_{xx}^3)+O(t^2). 
 \end{equation}
\normalsize
The coefficients in \cref{eq:hfunc} and \cref{eq:hder} obey the following recursive relation:
\small
\begin{align}
\begin{split}
\label{eq:rec}
&\beta_{1,0}=\beta_{2,0}=\beta_{3,0}=\beta_{4,0}=\beta_{5,0}=0\\
\beta_{1,m}=&\bigg[1+\sum_{i=0}^{m-1}\gamma_{m,i}\beta_{1,i}\bigg]\\
\beta_{2,m}=&\bigg[\frac{1}{2}+\sum_{i=0}^{m-1}\gamma_{m,i}\beta_{2,i}+\sum_{i=0}^{m-1}\gamma_{m,i}\beta_{4,i}\bigg]\\
\beta_{3,m}=&\bigg[\frac{2}{3}+\frac{1}{6}\bigg(\sum_{i=0}^{m-1}\gamma_{m,i}\beta_{1,i}\bigg)^3-\frac{1}{4}\bigg(\sum_{i=0}^{m-1}\gamma_{m,i}\beta_{1,i}\bigg)\bigg(\sum_{i=0}^{m-1}\gamma_{m,i}\beta_{1,i}^2\bigg)\\&+\frac{1}{12}\sum_{i=0}^{m-1}\gamma_{m,i}\beta_{1,i}^3+\sum_{i=0}^{m-1}\gamma_{m,i}\beta_{3,i}+\sum_{i=0}^{m-1}\gamma_{m,i}\beta_{5,i} \bigg]\\
\beta_{4,m}=&\bigg[1+\sum_{i=0}^{m-1}\gamma_{m,i}\beta_{4,i}\bigg] \\
\beta_{5,m}=&\bigg[2+\sum_{i=0}^{m-1}\gamma_{m,i}\beta_{5,i}\bigg]
\end{split}
\end{align}
\normalsize
Furthermore, the following conditions are necessary and
sufficient for second order accuracy of \cref{alg:mstage}:

\begin{align}
\label{eq:cons}
\begin{split}
&\frac{\beta_{2,M}}{\beta_{1,M}^2}=\frac{1}{2}\\
&\frac{\beta_{3,M}}{\beta_{1,M}^2}=1
\end{split}
\end{align}

\end{claim}
\begin{proof}
We will prove \cref{eq:hfunc} and \cref{eq:hder} by induction. For $h_1$, these equations are the two dimensional version of equations \cref{eq:vbar} and \cref{eq:vbarder} worked out in the previous section.\\
For the induction step, assume \cref{eq:hfunc} and \cref{eq:hder} up to $m-1$. We want to solve for $y$ such that $\bigg[G_t*\sum_{i=0}^{m-1} \gamma_{m,i} \bar{u}_m\bigg](0,y)=\frac{1}{2}$
Using \cref{eq:fconv} for two dimensions and the linearity of the  convolution we arrive at:
\small
\begin{align*}
\bigg[G_t&*\sum_{i=0}^{m-1} \gamma_{m,i} \bar{u}_i\bigg](0,y)\\=&\frac{1}{2}-\frac{y}{2\sqrt{\pi t}}+\frac{y^3}{24\sqrt{\pi }t^{3/2}}\\&+\frac{1}{2\sqrt{\pi t}}\sum_{i=0}^{m-1} \gamma_{m,i} \bigg[h_i +t \frac{d^2}{dx^2} h_i+\frac{t^2}{2} \frac{d^4}{dx^4} h_i-\frac{y^2}{4t}h_i-\frac{y^2}{4}\frac{d^2}{dx^2} h_i+\frac{y}{4t}h_i^2\\&+\frac{y}{2}h_i\frac{d^2}{dx^2}h_i+\frac{3ty}{4}(\frac{d^2}{dx^2}h_i)^2-\frac{h_i^3}{12t}-\frac{h_i^2}{4}\frac{d^2}{dx^2}h_i-\frac{3t}{4}h_i(\frac{d^2}{dx^2}h_i)^2-\frac{5t^2}{4}(\frac{d^2}{dx^2}h_i)^3\bigg]\\&+\text{h.o.t.} \\
=&\frac{1}{2}-\frac{y}{2\sqrt{\pi t}}+\frac{y^3}{24\sqrt{\pi }t^{3/2}}\\&+\frac{1}{2\sqrt{\pi t}}\sum_{i=0}^{m-1} \gamma_{m,i} \bigg[t \beta_{1,i} f_{xx}+t^2\beta_{2,i}f_{xxxx}-t^2\beta_{3,i}f_{xx}^3 +tf_{xx}+t^2\beta_{4,i}f_{xxxx}-t^2\beta_{5,i}f_{xx}^3\\&+\frac{t^2}{2}f_{xxxx}-\frac{y^2}{4}\beta_{1,i} f_{xx}-\frac{y^2}{4}f_{xx}+\frac{yt}{4}\beta_{1,i}^2f_{xx}^2+\frac{yt}{2}\beta_{1,i}f_{xx}^2\\&+\frac{3ty}{4}f_{xx}^2-\frac{\beta_{1,i}^3t^2}{12}f_{xx}^3-\frac{\beta_{1,i}^2t^2}{4}f_{xx}^3-\frac{3\beta_{1,i}t^2}{4}f_{xx}^3-\frac{5}{4}t^2f_{xx}^3\bigg] +\text{h.o.t.}
\end{align*}
\normalsize
Setting the previous equation equal to a half and solving for $y$ we have
\small
\begin{align}
\label{eq:solvey}
\begin{split}
y=&tf_{xx}\bigg[1+\sum_{i=0}^{m-1}\gamma_{m,i}\beta_{1,i}\bigg]+t^2f_{xxxx}\bigg[\frac{1}{2}+\sum_{i=0}^{m-1}\gamma_{m,i}\beta_{2,i}+\sum_{i=0}^{m-1}\gamma_{m,i}\beta_{4,i}\bigg]\\
&-t^2f_{xx}^3\bigg[\frac{2}{3}+\frac{1}{6}\bigg(\sum_{i=0}^{m-1}\gamma_{m,i}\beta_{1,i}\bigg)^3-\frac{1}{4}\bigg(\sum_{i=0}^{m-1}\gamma_{m,i}\beta_{1,i}\bigg)\bigg(\sum_{i=0}^{m-1}\gamma_{m,i}\beta_{1,i}^2\bigg)\\&+\frac{1}{12}\sum_{i=0}^{m-1}\gamma_{m,i}\beta_{1,i}^3+\sum_{i=0}^{m-1}\gamma_{m,i}\beta_{3,i}+\sum_{i=0}^{m-1}\gamma_{m,i}\beta_{5,i} \bigg]+\text{h.o.t.}
\end{split}
\end{align}
\normalsize
Similarly, using \cref{eq:der2} for two dimensions and the linearity of the  convolution we derive:
\small
\begin{align}
\label{eq:msder}
\begin{split}
\frac{d^2}{d x^2}& \bigg[G_t*\sum_{i=0}^{m-1} \gamma_{m,i} \bar{u}_i \bigg](x,y(x))|_{x=0}\\
=&-\frac{y_{xx}}{2\sqrt{\pi t}}+\frac{y^2y_{xx}}{8\sqrt{\pi}t^{3/2}}\\
&+\sum_{i=0}^{m-1} \gamma_{m,i}\bigg[\frac{1}{2\sqrt{\pi t}}\frac{d^2}{dx^2}h_i+\frac{\sqrt{t}}{2\sqrt{\pi}} \frac{d^4}{d x^4} h_i-\frac{y^2}{8\sqrt{\pi}t^{3/2}}(\frac{d^2}{dx^2}h_i)+\frac{y}{4\sqrt{\pi}t^{3/2}}h_i \bigg( \frac{d^2}{dx^2}h_i\bigg)\\&+\frac{3y}{4\sqrt{\pi t}}\bigg(\frac{d^2}{dx^2} h_i \bigg)^2 -\frac{h_i^2}{8\sqrt{\pi}t^{3/2}}\frac{d^2}{dx^2} h_i-\frac{3h_i}{4\sqrt{\pi t}}\bigg(\frac{d^2}{dx^2} h_i \bigg)^2\\&-\frac{15\sqrt{t}}{8\sqrt{\pi}}\bigg(\frac{d^2}{dx^2} h_{xx}\bigg)^3+\frac{y_{xx}h_i^2}{8\sqrt{\pi}t^{3/2}}+\frac{y_{xx}}{4\sqrt{\pi t}}h_i\bigg(\frac{d^2}{dx^2} h_i\bigg)+\frac{3\sqrt{t}}{8\sqrt{\pi}}t y_{xx}\bigg(\frac{d^2}{dx^2} h_i\bigg)^2\\&-\frac{y_{xx}yh_i}{4\sqrt{\pi t}}-\frac{y_{xx}y}{4\sqrt{\pi t}}\frac{d^2}{dx^2} h_i\bigg]+\text{h.o.t.}\\
\end{split}
\end{align}
Now substitute in \cref{eq:solvey}, \cref{eq:hfunc}, \cref{eq:hder} for $y$, $h_i$ and $\frac{d^2}{dx^2} h_i$ respectively in \cref{eq:msder}:
\small
\begin{align}
\label{eq:msderfinal}
\begin{split}
\frac{d^2}{d x^2}& \bigg[G_t*\sum_{i=0}^{m-1} \gamma_{m,i} \bar{u}_i \bigg](x,y(x))|_{x=0}\\
=&-\frac{y_{xx}}{2\sqrt{\pi t}}+\frac{1}{8\sqrt{\pi t}}\bigg(1+\sum_{i=0}^{m-1} \gamma_{m,i} \beta_{1,i}\bigg)^2y_{xx}f_{xx}^2+\frac{1}{2\sqrt{\pi t}}f_{xx}+\frac{\sqrt{t}}{2\sqrt{\pi}}\bigg(\sum_{i=0}^{m-1} \gamma_{m,i}\beta_{4,i}\bigg)f_{xx}^3\\
&+\frac{\sqrt{t}}{2\sqrt{\pi }}\bigg(\sum_{i=0}^{m-1} \gamma_{m,i}\beta_{5,i}\bigg)f_{xxxx}+\frac{\sqrt{t}}{2\sqrt{\pi}} f_{xxxx}-\frac{\sqrt{t}}{8\sqrt{\pi}}\bigg(1+\sum \gamma_{m,i} \beta_{1,i}\bigg)^2f_{xx}^3\\&+\frac{\sqrt{t}}{4\sqrt{\pi}}\bigg(1+\sum_i^{m-1} \gamma_{m,i} \beta_{1,i}\bigg)\bigg(\sum \gamma_{m,i} \beta_{1,i}\bigg)f_{xx}^3+\frac{3\sqrt{t}}{4\sqrt{\pi }}\bigg(1+\sum \gamma_{m,i} \beta_{1,i}\bigg)f_{xx}^3 \\&-\frac{\sqrt{t}}{8\sqrt{\pi}}\bigg(\sum_{i=0}^{m-1} \gamma_{m,i} \beta_{1,i}\bigg)^2f_{xx}^3-\frac{3\sqrt{t}}{4\sqrt{\pi}}\bigg(\sum_i^{m-1} \gamma_{m,i} \beta_{1,i}\bigg)f_{xx}^3-\frac{15\sqrt{t}}{8\sqrt{\pi}}f_{xx}^3 \\&+\frac{\sqrt{t}}{8\sqrt{\pi}}\bigg(\sum_{i=0}^{m-1} \gamma_{m,i} \beta_{1,i}\bigg)^2y_{xx}f_{xx}^2+\frac{\sqrt{t}}{4\sqrt{\pi }}\bigg(\sum_{i=0}^{m-1} \gamma_{m,i} \beta_{1,i}\bigg)y_{xx}f_{xx}^2+\frac{3\sqrt{t}}{8\sqrt{\pi}}y_{xx}f_{xx}^2\\&-\frac{\sqrt{t}}{4\sqrt{\pi}}\bigg(1+\sum_{i=0}^{m-1} \gamma_{m,i} \beta_{1,i}\bigg)\bigg(\sum_{i=0}^{m-1} \gamma_{m,i} \beta_{1,i}\bigg)y_{xx}f_{xx}^2-\frac{\sqrt{t}}{4\sqrt{\pi }}\bigg(1+\sum_{i=0}^{m-1} \gamma_{m,i} \beta_{1,i}\bigg)y_{xx}f_{xx}^2\\
&+\text{h.o.t.}
\end{split}
\end{align}
\normalsize
Setting \cref{eq:msderfinal} to zero and solving for $y_{xx}$ we find
\small
\begin{equation}
\label{eq:dery}
y_{xx}=f_{xx}+t\bigg[\bigg(1+\sum_{i=0}^{m-1}\gamma_{m,i}\beta_{4,i}\bigg)f_{xxxx}-\bigg(2+\sum_{i=0}^{m-1}\gamma_{m,i}\beta_{5,i}\bigg)f_{xx}^3\bigg]+O(t^2).
\end{equation}
\normalsize
Equations \cref{eq:solvey} and \cref{eq:dery} give the recursive relations \cref{eq:rec}.\\

The consistency equations \cref{eq:cons} follow by the change of variable $\tau=t \beta_{1,M}$ for $h_M(0)$  and matching the Taylor expansion for motion by mean curvature \cref{eq:true}. 
\end{proof}

\subsection{Unconditional Stability}

Next, we give conditions on the $\gamma$'s that preserve unconditional stability in any dimension. Specifically, for energy
\begin{equation}
\label{eq:energyn}
    E_t(\Sigma)=\int_{\mathbb{R}^n} (1-\mathbf{1}_{\Sigma})G_t*\mathbf{1}_{\Sigma} d\textbf{x}.
\end{equation}
our algorithm has the property $E_t(u_{n+1}) \leq E_t(u_n)$.
In \cite{variational2019}, the authors prove conditions for unconditional stability of the following class of linear $M$-stage algorithms:
\begin{equation}
\label{eq:ms}
    u_{n+1}=U_M=\argmin_u E(u)+\sum^{M-1}_{i=0}\frac{\gamma_{M,i}}{2k}||u-U_i||^2
\end{equation}
where the intermediate stages $U_m$, for $m \geq 1$, are given by
\begin{equation}
\label{eq:msintermediate}
U_m=\argmin_u E(u)+\sum^{m-1}_{i=0}\frac{\gamma_{m,i}}{2k}||u-U_i||^2.
\end{equation}
for some energy $E$, fixed time step $k$ and the stipulation $U_0=u_n$.
We state the stability conditions from that paper below, and show that \cref{alg:mstage} falls into the desired class:

\begin{theorem}(From \cite{variational2019})
\label{claim:ms}
Define the following auxiliary quantities in terms of the coefficients $\gamma_{m,i}$ of scheme \cref{eq:ms} and \cref{eq:msintermediate}:
\begin{align}
&\tilde{\gamma}_{m,i}=\gamma_{m,i}-\sum_{j=m+1}^M\tilde{\gamma}_{j,i}\frac{\tilde{S}_{j,m}}{\tilde{S}_{j,j}}\\
\label{eq:S}
&\tilde{S}_{j,m}=\sum_{i=0}^{m-1} \tilde{\gamma}_{j,i}
\end{align}
If $\tilde{S}_{m,m}>0$ for $m=1,\ldots,M$, then scheme \cref{eq:tdmsenergy} ensures that for every $n=0,1,2,\ldots$ we have $E(u_{n+1}) \leq E(u_n)$.
\end{theorem}

It is shown in \cite{eo} that each step of the original threshold dynamics, \cref{alg:mbo}, can be written as
\begin{equation*}
\Sigma_{n+1}=\argmin_{\Sigma} \int (1-\mathbf{1}_{\Sigma})G_t*\mathbf{1}_{\Sigma} + \int (\mathbf{1}_{\Sigma}-\mathbf{1}_{\Sigma_n})G_t*(\mathbf{1}_{\Sigma}-\mathbf{1}_{\Sigma_n}) \, d\mathbf{x}.
\end{equation*}
Similarly observe that \cref{eq:tdms} can be written as
\begin{equation}
\label{eq:tdmsenergy}
\bar{\Sigma}_m=\argmin_{\Sigma} \int (1-\mathbf{1}_{\Sigma})G_t*\mathbf{1}_{\Sigma} +\sum_{i=0}^{m-1} \gamma_{m,i} \int (\mathbf{1}_{\Sigma}-\mathbf{1}_{\bar{\Sigma}_i})G_t*(\mathbf{1}_{\Sigma}-\mathbf{1}_{\bar{\Sigma}_i}) \, d\mathbf{x}.
\end{equation}
Moreover, as noted in \cite{eo}, since $\widehat{G}_t>0$,
\[\int u G_t * u \geq 0 \, d\mathbf{x}\] 
with equality holding if and only if $u=0$.
Thus, \cref{alg:mstage} is of type \cref{eq:ms} and \cref{eq:msintermediate} with energy \cref{eq:energyn} and inner product $\langle u , v \rangle = \int u G_t * v \, d\mathbf{x}$, so that conditions for unconditional stability identified in \cite{variational2019} apply.\\

In the next section, we give examples of $\gamma$'s that satisfy the consistency equations (\cref{claim:cons}) as well as the hypothesis of \cref{claim:ms} concurrently. 

\subsection{A Second Order Unconditionally Stable Example}
\label{sec:example}
In this section we present a set of $\gamma$'s such that the second order consistency equations \cref{eq:cons} and hypothesis of unconditional stability \cref{claim:ms} are satisfied simultaneously.  We found the $\gamma$'s  numerically and then sought a nearby algebraic solution to the consistency equations that still satisfied the conditions in \cref{claim:ms}. It can be shown that there is no unconditionally energy stable second order two- or three-stage method of type \cref{eq:tdms} for threshold dynamics. However, a four-step method exists with the following $\gamma$'s:
\begin{equation}
\label{eq:gammas}
    \gamma=\left(
\begin{array}{cccc}
 \gamma_{1,0} & 0 & 0 & 0 \\
 \gamma_{2,0}& \gamma_{2,1}& 0 &0 \\
\gamma_{3,0} & \gamma_{3,1} & \gamma_{3,2} &0 \\
  \gamma_{4,0} & \gamma_{4,1} & \gamma_{4,2} &\gamma_{4,3}\\
\end{array}\right) \\ \approx \left(
\begin{array}{cccc}
 1 & 0 & 0 & 0 \\
 -0.25 & 1.25 & 0 &0 \\
 0.83 & -0.67 & 0.83 &0 \\
  -0.73 & 0.5 & -0.73 &1.96\\
\end{array}\right).\\
\end{equation}
The exact values are in the appendix (\cref{sec:gamma2}); they are all algebraic numbers but with long representations. This choice of $\gamma$'s that endows \cref{alg:mstage} with unconditional stability and second order accuracy is not unique. In fact, one could find other choices of $\gamma$ that preserve additional desirable properties. We summarize our results in the following theorem:
\begin{theorem}
\cref{alg:mstage} with coefficients \cref{eq:gammas} is unconditionally energy stable and first order accurate in any dimension. Moreover, it is second order accurate in two dimensions.
\end{theorem} 
\subsection{Consistency In Higher Dimensions}
Unfortunately, there does not exist a choice of $\gamma$'s such that \cref{alg:mstage} is second order in three dimensions. Using the Gaussian (or even a linear combination of Gaussian kernels) at every step in our multistage algorithm leads to consistency equations incompatible with second order mean curvature motion \cref{eq:mctaylor}. Of course, one can choose a different kernel, denoted by $K_t$, at each stage. The consistency equations will be different. On the other hand, if the kernel has positive Fourier transform, $\widehat{K}_t\geq 0$, \cref{claim:ms} will hold for energy 
\[
E_t^*(\Sigma)=\int (1-\mathbf{1}_{\Sigma})K_t*\mathbf{1}_{\Sigma} d\mathbf{x}.
\]
As an additional consideration, the energy induced by the kernel should have property 
\begin{equation}
\label{eq:gammacon}
    \lim_{t \to 0^+} E^*_t(\Sigma) \xrightarrow{\Gamma} c\text{Per}(\Sigma)
\end{equation} for some constant $c$. (The Gaussian kernel has this property~\cite{m77}.) It is left to future work to find a scheme with kernel and $\gamma$'s such that unconditional stability, positive Fourier transform, and property \cref{eq:gammacon} hold, and that, furthermore, is second order in three (and higher) dimensions. 

\section{Numerical Tests}
\label{sec:numerical}
In this section, we present highly accurate convergence tests for the two new algorithms: \cref{alg:twostep} and \cref{alg:mstage}.
It is well known that naive implementations of threshold dynamics on uniform grids introduces large spatial discretization errors due to sampling characteristic functions.
In fact, if the time step size is sufficiently small compared to the spatial discretization, interfaces can even get stuck.
A very effective cure to this phenomenon is the adaptive implementation of Ruuth \cite{ruuth1998efficient}, which entails fast convolution on non-uniform grids.
In practice, we recommend that the high order in time schemes introduced in this paper be used in conjunction with such a spatial implementation.

Below, our focus is on numerically verifying the improvement in the convergence rate in time of the new threshold dynamics schemes on a few smooth interfaces.
To ensure spatial errors are negligible in our experiments, for most of our experiments below, we simply represent the interfaces as graphs of functions.
This is, of course, not meant as a practical implementation, but just as a very accurate and efficient way to minimize spatial errors -- it allows us to reach very high spatial resolutions and accuracies -- so that we can focus on time discretization errors.
Section \ref{sec:numerical.1} explains the details of the implementation, and Section \ref{sec:numerical.2} presents the results.
The latter also contains an experiment with the practical implementation of the new schemes (albeit on uniform grids) to verify that no substantial qualitative difference is observed in the handling of topological changes vs. the original threshold dynamics algorithm.
\subsection{Highly Accurate Threshold Dynamics For Graphs}
\label{sec:numerical.1}
This section explains how \cref{alg:twostep} and \cref{alg:mstage} have been implemented in the special case that the interface is given as the graph of a function for the purpose of numerical convergence tests given in Section \ref{sec:numerical.2}.
We repeat that we are not advocating this implementation as a practical method -- it is just for numerical convergence tests on smooth interfaces -- but refer to Ruuth's adaptive spectral implementation \cite{ruuth1998efficient} to be used in conjunction with our new algorithms.
\\

Let $\Sigma=\{(\mathbf{x},x_n)| x_n \leq f(\mathbf{x})\}$ for some function $f(\mathbf{x})$.
Recall the definition of the Gaussian kernel $G_{\delta t}(\mathbf{x})=\frac{1}{(4 \pi \delta t)^{d/2}}e^{-\frac{|\mathbf{x}|^2}{4\delta t}}$. Then the convolution is
\begin{align}
\label{eq:approx}
\begin{split}
G_{\delta t} * \mathbbm{1}_{\Sigma}(\mathbf{x},x_n)&=\int_{\mathbb{R}^{d-1}}G_{\delta t}(\mathbf{x}-\bar{\mathbf{x}}) \int_{-\infty}^{f(\bar{\mathbf{x}})} \frac{1}{\sqrt{4 \pi \delta t}}e^{-\frac{(x_n-\bar{x}_n)^2}{4\delta t}} d\bar{\mathbf{x}} d\bar{x}_n\\
&=\frac{1}{2}
+\frac{1}{2}\int_{\mathbb{R}^{d-1}}G_{\delta t}(\mathbf{x}-\bar{\mathbf{x}}) \erf\bigg(\frac{f(\bar{\mathbf{x}})-x_n}{2\sqrt{t}}\bigg) d\bar{\mathbf{x}}
\end{split}
\end{align}
where $\erf(\cdot)$ is the error function. The latter integral is calculated numerically (e.g. Gaussian quadrature) inside a region $\Omega \subset \mathbb{R}^{n-1}$ such that $\int_{\Omega^c}G_{\delta t}(\mathbf{x}-\bar{\mathbf{x}})d\bar{\mathbf{x}} <\epsilon$ for some preset tolerance $\epsilon$.

With this tool in hand, we can implement the original threshold dynamics (\cref{alg:mbo}) by tracking the interface at fixed points $\mathbf{x}$ throughout the evolution. We describe this in \cref{alg:threspar}. The two version of threshold dynamics we proposed in this paper,\cref{alg:twostep} and \cref{alg:mstage}, are implemented similarly. 

\begin{algorithm}
  \caption{Highly Accurate Threshold Dynamics For Graphs
    \label{alg:threspar}}
  \begin{algorithmic}[1]
  \State Fix total evolution time $T$, time step size $\delta t$, and points $\{\mathbf{x}^i\}_{i=1}^N \in D \subset \mathbb{R}^{d-1}$.
  \State Set $\Sigma^0=\{(\mathbf{x},x_n)| x_n \leq f(\mathbf{x},0)\}$ and $K=T/\delta t$.
      \For{$k \gets 1 \textrm{ to } K$}
      	  \State For each $\mathbf{x}_i$ find $y^i$ such that \[\frac{1}{2}=G_{\delta t} * \mathbbm{1}_{\Sigma^{k-1}}(\mathbf{x}^i,y^i)\] using e.g. the secant method,
      	  estimating the right hand side using \cref{eq:approx}.
          \State Set $f(\mathbf{x}^i,k\delta t)=y^i$ and  $\Sigma^k=\{(\mathbf{x},x_n)| x_n \leq f(\mathbf{x},k\delta t)\}$.
			
      \EndFor
  \end{algorithmic}{}
\end{algorithm}

\subsection{Numerical Results}
\label{sec:numerical.2}
In this section, we will test \cref{alg:twostep} and \cref{alg:mstage} on a couple of two phase mean curvature motion evolution problems to demonstrate their accuracy.\\

In $\mathbb{R}^2$ we run our algorithms on the `Grim Reaper Wave.' The interface is given by $f(x,t)= \arcsinh(\exp(-t)\cos(x))$. We run the evolution with initial data $f(x,0)=\arcsinh(\cos(x))$ to $T=1$ on domain $x=[0,\pi]$ with homogeneous Neumann boundary conditions. We track $4000$ points and report the $L^2$ error between final interface and $f(x,1)$. We include the errors for the original threshold dynamics in \cref{tab:TD22d} for comparison. The error for \cref{alg:twostep} is reported in \cref{tab:TD22dms} and the error for \cref{alg:mstage}, with $\gamma$'s given in \cref{eq:gammas}, is tabulated in \cref{tab:TD22dk}.

\begin{table}[h]
\begin{center}
\begin{tabular}{|c|c|c|c|c|c|c|}
\hline
Number of& & & & & & \\time steps&8&16 &32&64&128&256\\
\hline
$L^2$ error&5.4e-03&2.6e-03&1.3e-03&6.7e-04&3.3e-04&1.7e-04\\
\hline
Order&-&1.0&1.0&1.0&1.0&1.0\\
\hline
\end{tabular}
\caption{Ordinary Threshold Dynamics, \cref{alg:mbo}, on the `Grim Reaper Wave.'}
\label{tab:TD22d}
\end{center}
\end{table}

\begin{table}[h]
\begin{center}
\begin{tabular}{|c|c|c|c|c|c|c|}
\hline
Number of& & & & & & \\time steps&8&16 &32&64&128&256\\
\hline
$L^2$ error&1.9e-04& 3.4e-05 & 8.4e-06 &2.1e-06&5.1e-07& 1.3e-07\\
\hline
Order&-&2.5& 2.0&2.0& 2.0&  2.0\\
\hline
\end{tabular}
\caption{\cref{alg:twostep} on the `Grim Reaper Wave.'}
\label{tab:TD22dms}
\end{center}
\end{table}

\begin{table}[h]
\begin{center}
\begin{tabular}{|c|c|c|c|c|c|c|}
\hline
Number of& & & & & & \\time steps&8&16 &32&64&128&256\\
\hline
$L^2$ error&6.6e-05& 1.6e-05 & 4.1e-06 &1.0e-06&  2.6e-07&   6.3e-08\\
\hline
Order&-&2.0& 2.0&2.0& 2.0&  2.0\\
\hline
\end{tabular}
\caption{\cref{alg:mstage} with $\gamma$'s given in \cref{eq:gammas} on the `Grim Reaper Wave.'}
\label{tab:TD22dk}
\end{center}
\end{table}

In three dimensions, we run mean curvature evolution tests with initial interface $f(x,y,0)=\cos(\pi y)\cos(\pi x)+\frac{1}{2}\cos(\pi y)$ to $T=1/10$ on $(x,y) \in [-1,1]\times[-1,1]$ with periodic boundary conditions. We generate the `true' solution by forward Euler discretization of the PDE for mean curvature motion \cref{eq:mcmpde} using very small time steps and a very high spatial resolution. We tabulate the  $L^2$ error between the `true' interface and the output of \cref{alg:twostep} in \cref{tab:TD23d}.

\begin{figure}[h]
\begin{center}
\includegraphics[width=.32\textwidth]{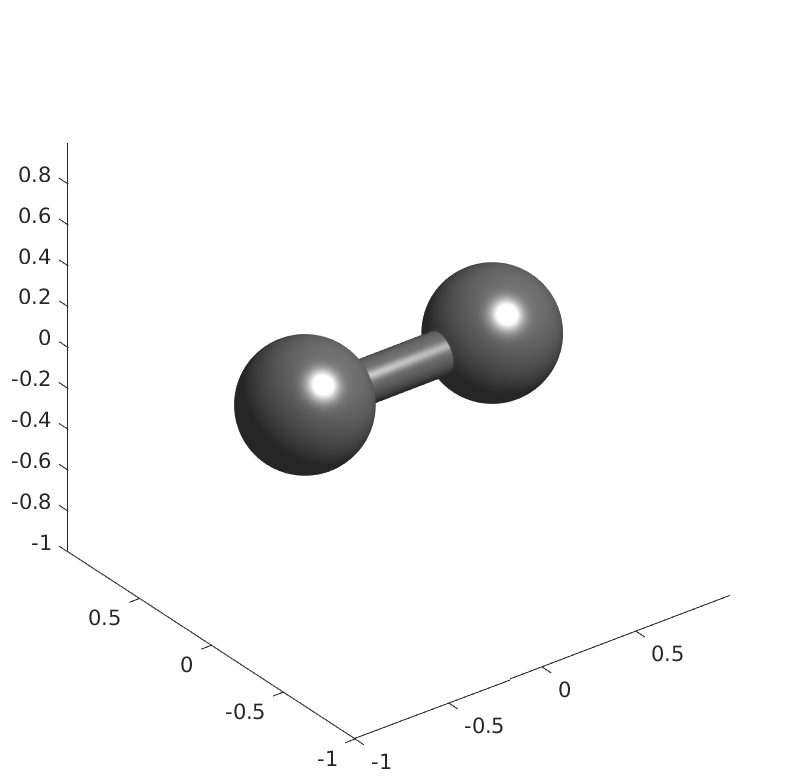}
\newline
\includegraphics[width=.32\textwidth]{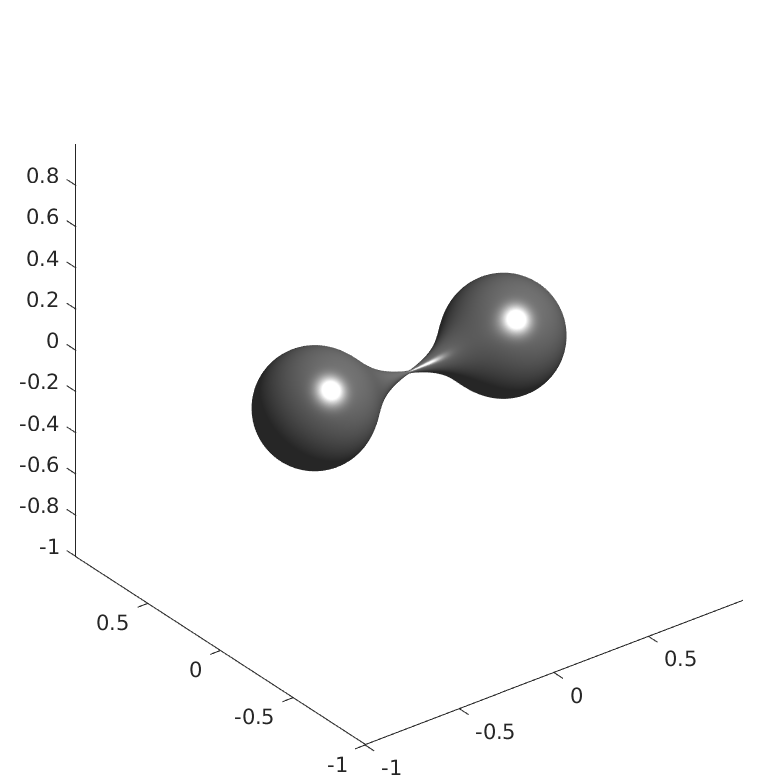}
\includegraphics[width=.32\textwidth]{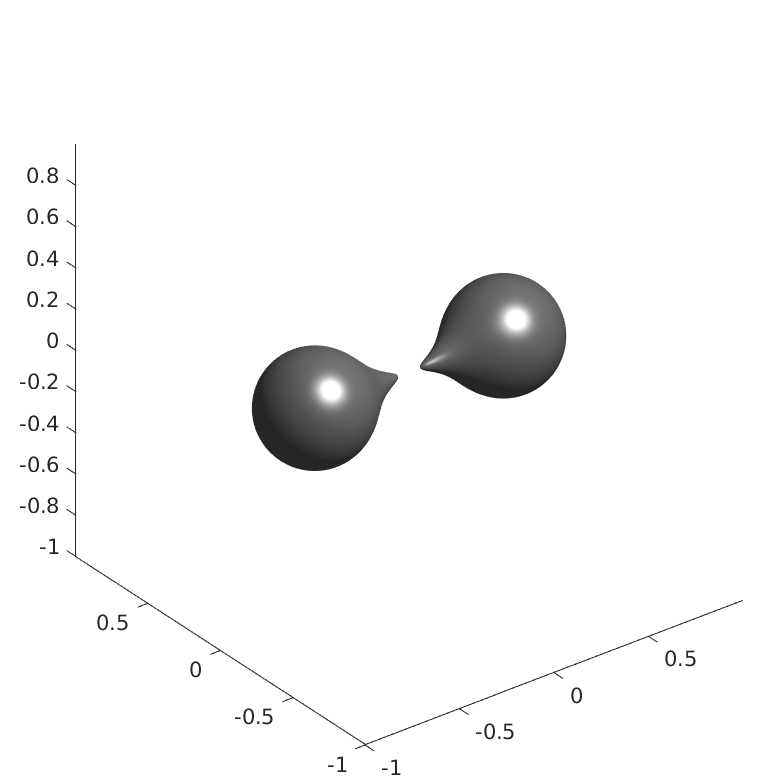}
\includegraphics[width=.32\textwidth]{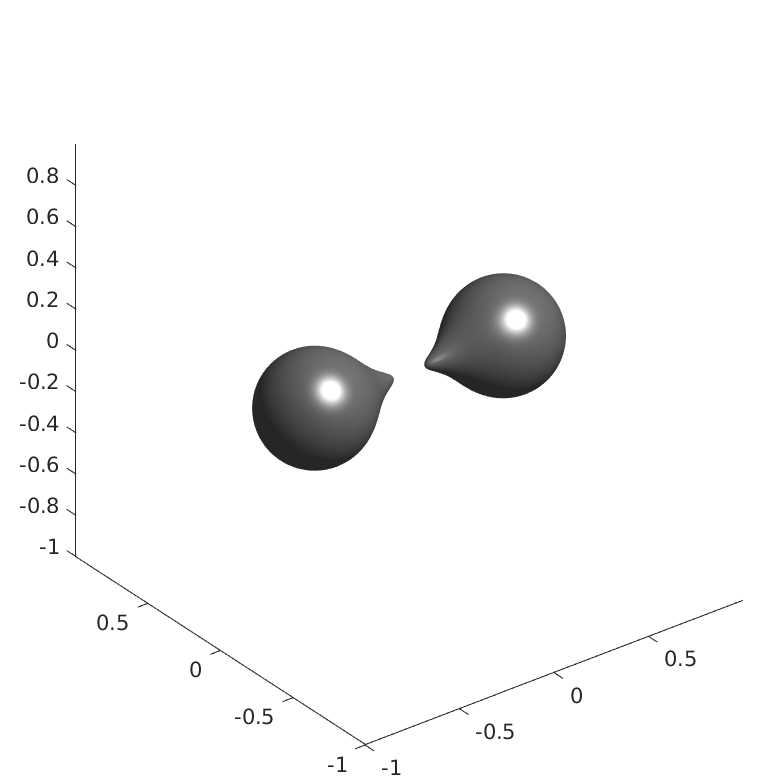}
\includegraphics[width=.32\textwidth]{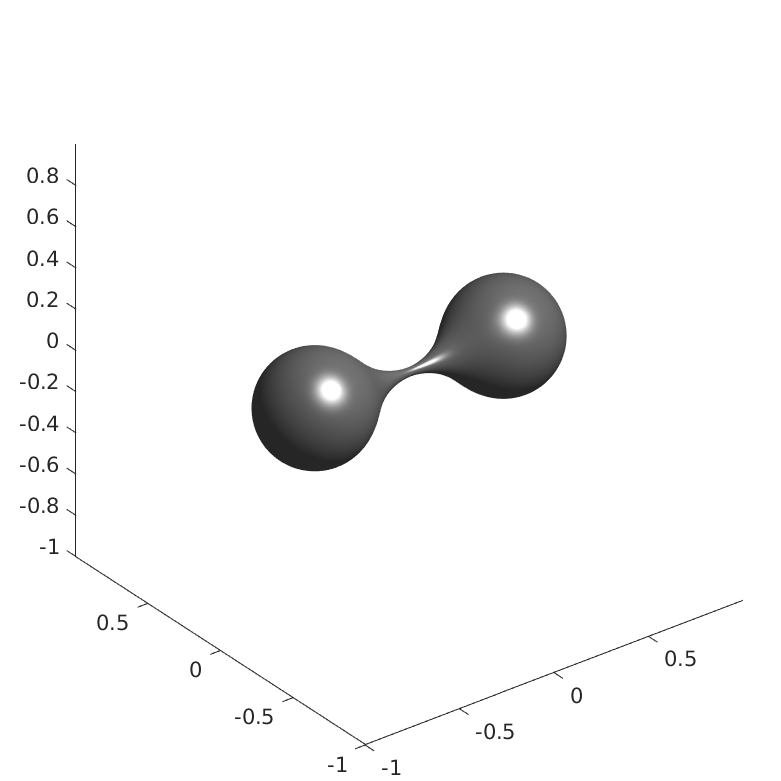}
\includegraphics[width=.32\textwidth]{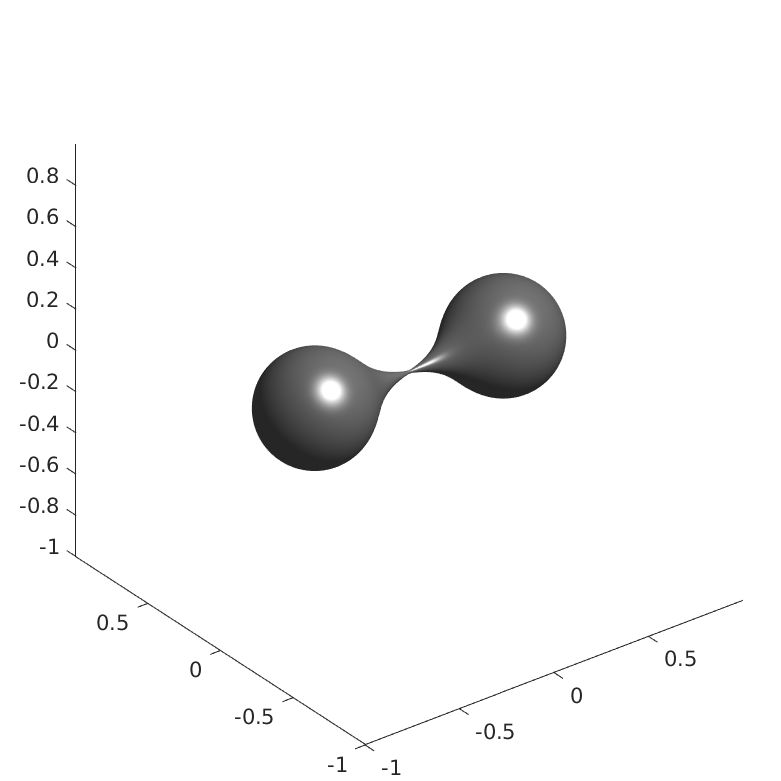}
\includegraphics[width=.32\textwidth]{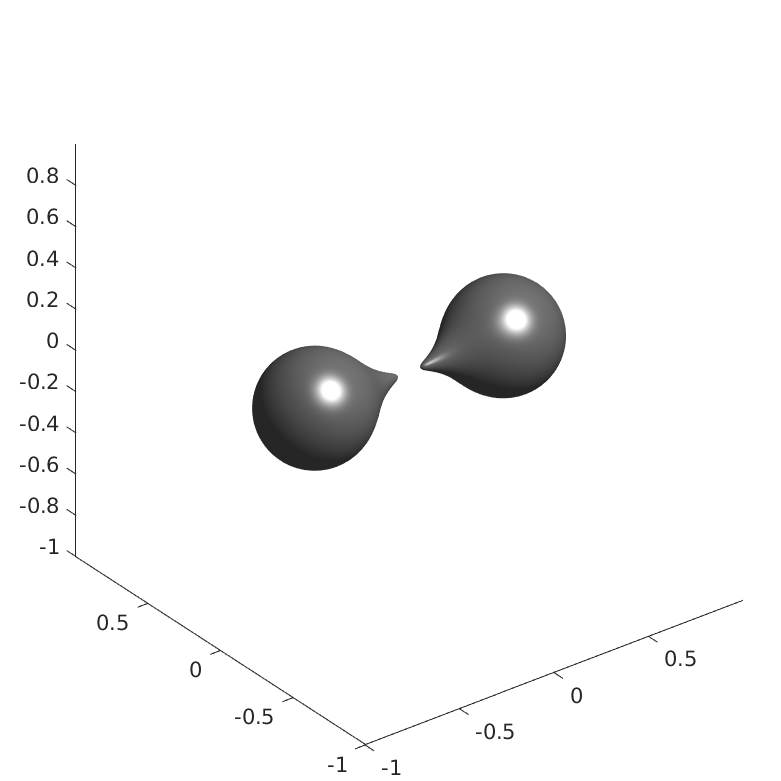}
\caption{\footnotesize The original threshold dynamics \cref{alg:mbo} and \cref{alg:twostep} evolving a dumbbell by mean curvature motion. Top: the initial dumbbell. Center: \cref{alg:mbo} before and after the topological change. Bottom: \cref{alg:twostep} before and after the topological change at the same time values as the center row.}
\label{fig:dumbbell}
\end{center}
\end{figure}

\begin{table}[h]
\begin{center}
\begin{tabular}{|c|c|c|c|c|c|}
\hline
Number of& & & & & \\time steps&4&8&16 &32&64\\
\hline
L2&3.3e-03&6.5e-04&1.3e-04&2.7e-05&6.3e-06\\
\hline
Order&-&2.4&2.4&2.2&2.1\\
\hline
\end{tabular}
\caption{\cref{alg:twostep} on an interface in $\mathbb{R}^3$}
\label{tab:TD23d}
\end{center}
\end{table}
As an additional test in three dimensions, we use \cref{alg:twostep} to evolve a dumbbell by mean curvature motion using the practical implementation of \cref{alg:twostep}, albeit on a uniform grid. The system goes through a topological change where the connected set breaks off into two components. In \cref{fig:dumbbell}, we show the original threshold dynamics \cref{alg:mbo} and \cref{alg:twostep} at time steps near the pinch off. There is a slight difference between the two algorithms at the temporal and spatial resolutions we have chosen. Presumably, this difference will shrink as we further refined our temporal and spatial resolution. Regardless, \cref{alg:twostep} behaves reasonably during the pinch off.

\section{Conclusion}
We have presented second order threshold dynamic methods for simulating two phase mean curvature flow. Unlike the previous method, \cref{alg:ruuth}, our methods represent the phases via a binary function at every step.  Additionally, we present an unconditionally stable method in two dimensions. We have demonstrated the methods and their advertised accuracy on examples in two and three dimensions.

Finding a method that is second order in three and higher dimensions and unconditionally stable in the case of two phases is a topic of future investigation.

\section{Appendix}
\subsection{Taylor Expansion of Characteristic Function with a Gaussian Kernel}
\label{sec:taylor}
In this section of the appendix, we work out the Taylor expansion of the convolution of a Gaussian kernel with a characteristic function. A simpler version of the following calculation is worked out in two dimensions by Ruuth \cite{Ruuth_1996} and up to first order in arbitrary dimensions by Grzhibovskis \& Heinz \cite{gh}.  First, let us introduce the following notation for the 1D Gaussian for convenience,
\[
g_t(x)=\frac{1}{2 \sqrt{\pi t}}\exp\bigg[-\frac{x^2}{4t}\bigg] \text{ and let }G_t(x,y,z)=g_t(x)g_t(y)g_t(z).
\]
Now take a function $h(x,y)$ with the following properties:
\small
\begin{align}
\label{eq:hcond}
\begin{split}
& h(0,0)=\mathcal{O}(t) \\
& h_x(0,0)=\mathcal{O}(t) \\
& h_y(0,0)=\mathcal{O}(t)
\end{split}
\end{align}
\normalsize
The papers mentioned above (\cite{gh},\cite{Ruuth_1996}) use the assumption that \begin{equation}
    \label{eq:simpleh}
    h(0,0)= h_x(0,0)=h_y(0,0)=0.
\end{equation} These simpler assumptions are sufficient for their calculation of the Taylor expansion for a single step of threshold dynamics. However, after the first stage, the interface may no longer satisfy \cref{eq:simpleh}. So we required the more general conditions given by \cref{eq:hcond} for Taylor expansions of the interface after the first stage.\\

Now let $\Sigma=\{(x,y,z):z\leq h(x,y)\}$. The goal is to see how $\bigg[G_t*\mathbbm{1}_{\Sigma}\bigg](x,y,z)$ behaves along the $z$-axis near the origin. 
First, we simplify $\bigg[G_t*\mathbbm{1}_{\Sigma}\bigg](x,y,z)$:
\small
\begin{align}
\label{eq:conv}
\begin{split}
\bigg[G_t&*\mathbbm{1}_{\Sigma}\bigg](x,y,z)\\ =&\int_{\mathbb{R}^2} g_t(y-\tilde{y}) g_t(x-\tilde{x}) \int_{-\infty}^{\infty} g_t(z-\tilde{z})\mathbbm{1}_{\Sigma}(\tilde{x},\tilde{y},\tilde{z}) d\tilde{z} d\tilde{x}d\tilde{y}\\
=&\int_{\mathbb{R}^2} g_t(y-\tilde{y}) g_t(x-\tilde{x}) \int_{-\infty}^{h(\tilde{x},\tilde{y})} g_t(z-\tilde{z})  d\tilde{z} d\tilde{x}d\tilde{y}\\
=&\int_{\mathbb{R}^{2}} g_t(y-\tilde{y}) g_t(x-\tilde{x}) \int_{-\infty}^{0} g_t(z-\tilde{z}) d\tilde{z} d\tilde{x}d\tilde{y}\\&+\int_{\mathbb{R}^{2}}  g_t(y-\tilde{y}) g_t(x-\tilde{x}) \int_{0}^{h(\tilde{x},\tilde{y})} g_t(z-\tilde{z}) d\tilde{z} d\tilde{x}d\tilde{y}\\
=&\frac{1}{2}-\frac{z}{2\sqrt{\pi t}}+\frac{z^3}{24\sqrt{\pi}t^{3/2}}+\int_{\mathbb{R}^{2}} g_t(x-\tilde{x})g_t(y-\tilde{y}) \int_{0}^{h(\tilde{x},\tilde{y})} g_t(z-\tilde{z}) d\tilde{z} d\tilde{x}d\tilde{y}\\&+\text{h.o.t.}
\end{split}
\end{align}
\normalsize
Setting $x=y=0$, we will now simplify the last term of \cref{eq:conv},
\small
\begin{equation}
    \label{eq:lastterm}
    \int_{\mathbb{R}^{2}} g_t(x-\tilde{x})g_t(y-\tilde{y}) \int_{0}^{h(\tilde{x},\tilde{y})} g_t(z-\tilde{z}) d\tilde{z} d\tilde{x}d\tilde{y}.
\end{equation}
\normalsize
First note that, near $z=0$,
\small
\begin{equation}
\label{eq:intapprox}
\int_{0}^{h(\tilde{x},\tilde{y})} g_t(z-\tilde{z}) d\tilde{z} =\frac{1}{2\sqrt{\pi t}}\int_{0}^{h(\tilde{x},\tilde{y})} 1-\frac{(z-\tilde{z})^2}{4t} d\tilde{z}+\text{h.o.t.}
\end{equation}
\normalsize
Substituting approximation \cref{eq:intapprox} into \cref{eq:lastterm} and integrating,
\small
\begin{align}
\label{eq:intaprox}
\begin{split}
&\frac{1}{2\sqrt{\pi t}}\int_{\mathbb{R}^{2}} g_t(\tilde{x})g_t(\tilde{y}) \int_{0}^{h(\tilde{x},\tilde{y})} g_t(z-\tilde{z}) d\tilde{z} d\tilde{x}d\tilde{y}\\
=&\frac{1}{2\sqrt{\pi t}}\int_{\mathbb{R}^{2}} g_t(\tilde{x}) g_t(\tilde{y}) \int_{0}^{h(\tilde{x},\tilde{y})} 1-\frac{(z-\tilde{z})^2}{4t} d\tilde{z}d\tilde{x}d\tilde{y}+\text{h.o.t.}\\
=&\frac{1}{2\sqrt{\pi t}}\int_{\mathbb{R}^{2}} g_t(\tilde{x}) g_t(\tilde{y}) \bigg[h(\tilde{x},\tilde{y})+\frac{-3h(\tilde{x},\tilde{y})z^2+3(h(\tilde{x},\tilde{y}))^2z-(h(\tilde{x},\tilde{y}))^3}{12t}\bigg]  d\tilde{z}d\tilde{x}d\tilde{y} +\text{h.o.t.}
\end{split}
\end{align}
\normalsize

Denote the Taylor expansion of $h(\tilde{x},\tilde{y})$ around $(0,0)$ as $P[h](\tilde{x},\tilde{y})$:
\small
\begin{multline}
\label{eq:taylorh}
P[h](\tilde{x},\tilde{y})=h+\tilde{x}h_x+\tilde{y}h_y+\frac{\tilde{x}^2}{2}h_{xx}+\tilde{x}\tilde{y}h_{xy}+\frac{\tilde{y}^2}{2}h_{yy}+\frac{\tilde{x}^3}{6}h_{xxx}+\frac{\tilde{x}^2\tilde{y}}{2}h_{xxy}+\frac{\tilde{x}\tilde{y}^2}{2}h_{xyy}\\+\frac{\tilde{y}^3}{2}h_{yyy}+\frac{\tilde{x}^4}{24}h_{xxxx}+\frac{\tilde{x}^3\tilde{y}}{6}h_{xxxy}+\frac{\tilde{x}^2\tilde{y}^2}{4}h_{xxyy}+\frac{\tilde{x}\tilde{y}^3}{6}h_{xyyy}+\frac{\tilde{y}^4}{24}h_{yyyy}+\text{h.o.t.}
\end{multline}
\normalsize
Where $h=h(0,0)$, $h_x=h_x(0,0)$, etc. in order to simplify notation. Now substitute the Taylor expansion of $h(\tilde{x},\tilde{y})$ about $(0,0)$ into \cref{eq:intaprox}:
\small
\begin{multline}
\label{eq:expanded}
  \frac{1}{2\sqrt{\pi t}}\int_{\mathbb{R}^{2}} g_t(\tilde{x})g_t(\tilde{y}) \int_{0}^{h(\tilde{x},\tilde{y})} g_t(z-\tilde{z}) d\tilde{z} d\tilde{x}d\tilde{y}  =\\
  \frac{1}{2\sqrt{\pi t}}\int_{\mathbb{R}^{2}} g_t(\tilde{x}) g_t(\tilde{y}) \Bigg[P[h](\tilde{x},\tilde{y}) -\frac{1}{4t}P[h](\tilde{x},\tilde{y}) z^2 +\frac{1}{4t}\big(P[h](\tilde{x},\tilde{y})\big)^2z \\-\frac{1}{12t}\big(P[h](\tilde{x},\tilde{y})\big)^3\Bigg]  d\tilde{x}d\tilde{y}+\text{h.o.t.}
\end{multline}
\normalsize
Now we can integrate \cref{eq:expanded}. For a non-negative integer $n$ we have:
\small
\begin{equation} 
\label{eq:expmom}
\int_{-\infty}^{\infty} x^ng_t(x) dx=\begin{cases} 
      \frac{(2n)!}{n!}t^{n/2} & x\text{ for $n$ even} \\
      0 & \text{ for $n$ odd} 
   \end{cases}
\end{equation}
\normalsize
Using \cref{eq:expmom}, we simplify \cref{eq:expanded}
\small
\begin{align}
\begin{split}
\label{eq:exp}
\frac{1}{2\sqrt{\pi t}}\int_{\mathbb{R}^{2}} &g_t(\tilde{x})g_t(\tilde{y}) \int_{0}^{h(\tilde{x},\tilde{y})} g_t(z-\tilde{z}) d\tilde{z} d\tilde{x}d\tilde{y}=\\
\phantom{=}&\frac{h}{2\sqrt{\pi t}}+\frac{\sqrt{t}}{2\sqrt{\pi}}(h_{xx}+h_{yy})\\&+\frac{t^{3/2}}{4\sqrt{\pi}}(h_{xxxx}+2h_{xxyy}+h_{yyyy})
    -\frac{z^2}{8\sqrt{\pi}t^{3/2}}h-\frac{z^2}{8\sqrt{\pi t}}(h_{xx}+h_{yy})+\frac{z}{8\sqrt{\pi }t^{3/2}}h^2\\&+\frac{z}{4\sqrt{\pi t}}h(h_{xx}+h_{yy})+\frac{z\sqrt{t}}{2\sqrt{\pi}} (\frac{3}{4}h_{xx}^2+\frac{3}{4}h_{yy}^2+\frac{1}{2}h_{xx}h_{yy}+h_{xy}^2)-\frac{h^3}{24\sqrt{\pi}t^{3/2}} \\&-\frac{h^2}{8\sqrt{\pi t}}(h_{xx}+h_{yy})-\frac{\sqrt{t}}{2\sqrt{\pi}}h(\frac{3}{4}h_{xx}^2+\frac{3}{4}h_{yy}^2+\frac{1}{2}h_{xx}h_{yy}+h_{xy}^2)\\&-\frac{t^{3/2}}{2\sqrt{\pi}}(\frac{5}{4} h_{xx}^3+\frac{5}{4} h_{yy}^3+\frac{3}{4}h_{xx}h_{yy}^2+\frac{3}{4}h_{xx}^2h_{yy}+3h_{xx}h_{xy}^2+3h_{yy}h_{xy}^2)+\text{h.o.t.}
    \end{split}
\end{align}
\normalsize
Now substituting \cref{eq:exp} for the last term of \cref{eq:conv} we arrive at the expansion of $\bigg[G_t*\mathbbm{1}_{\Sigma}\bigg](0,0,z)$ near $z=0$:
\small
\begin{align}
\label{eq:fconv}
\begin{split}
\bigg[G_t&*\mathbbm{1}_{\Sigma}\bigg](0,0,z)=\\
&\frac{1}{2}-\frac{z}{2\sqrt{\pi t}}+\frac{z^3}{24\sqrt{\pi }t^{3/2}}+\frac{h}{2\sqrt{\pi t}}+\frac{\sqrt{t}}{2\sqrt{\pi}}(h_{xx}+h_{yy})\\&+\frac{t^{3/2}}{4\sqrt{\pi}}(h_{xxxx}+2h_{xxyy}+h_{yyyy})
    -\frac{z^2}{8\sqrt{\pi}t^{3/2}}h-\frac{z^2}{8\sqrt{\pi t}}(h_{xx}+h_{yy})+\frac{z}{8\sqrt{\pi }t^{3/2}}h^2\\&+\frac{z}{4\sqrt{\pi t}}h(h_{xx}+h_{yy})+\frac{z\sqrt{t}}{2\sqrt{\pi}} (\frac{3}{4}h_{xx}^2+\frac{3}{4}h_{yy}^2+\frac{1}{2}h_{xx}h_{yy}+h_{xy}^2)-\frac{h^3}{24\sqrt{\pi}t^{3/2}} \\&-\frac{h^2}{8\sqrt{\pi t}}(h_{xx}+h_{yy})-\frac{\sqrt{t}}{2\sqrt{\pi}}h(\frac{3}{4}h_{xx}^2+\frac{3}{4}h_{yy}^2+\frac{1}{2}h_{xx}h_{yy}+h_{xy}^2)\\&-\frac{t^{3/2}}{2\sqrt{\pi}}(\frac{5}{4} h_{xx}^3+\frac{5}{4} h_{yy}^3+\frac{3}{4}h_{xx}h_{yy}^2+\frac{3}{4}h_{xx}^2h_{yy}+3h_{xx}h_{xy}^2+3h_{yy}h_{xy}^2)+\text{h.o.t}
\end{split}
\end{align}
\normalsize

In this paper, we use the previous calculation to find the location of an interface along the $z$-axis after thresholding, i.e. finding $z$ such that  $\bigg[G_t*\mathbbm{1}_{\Sigma}\bigg](0,0,z)=\frac{1}{2}$. We also need to find how the derivatives of our interface, given by $z_x(0,0)$, $z_y(0,0)$, $z_{xx}(0,0)$ etc., relate to the derivatives of the  original interface given by $h(x,y)$. 
The derivatives of $z$ match the derivatives of $h$ to order $t$: 
\begin{equation}
\label{eq:blakentder}
\frac{\partial^{m+n}}{\partial x^n y^m} z(x,y)|_{(x,y)=(0,0)}=\frac{\partial^{m+n}}{\partial x^n y^m} h(x,y)|_{(x,y)=(0,0)}+O(t).
\end{equation}
Note that $z_x(0,0)$ and $z_y(0,0)$ are of $\mathcal{O}(t)$ for $h$ having properties \cref{eq:hcond}. For our calculations, we also need to find $z_{xx}(0,0)$ and $z_{yy}(0,0)$ to $\mathcal{O}(t^2)$.  
So we also include the calculation of $\frac{\partial^2}{\partial x^2} \bigg[G_t*\mathbbm{1}_{\Sigma}\bigg](x,y,z(x,y))|_{(x,y)=(0,0)}$:
\small
\begin{align}
\label{eq:conder}
\begin{split}
& \frac{\partial^2}{\partial x^2} \bigg[G_t*\mathbbm{1}_{\Sigma}\bigg](x,y,z(x,y))|_{(x,y)=(0,0)}\\
=& -\frac{z_{xx}}{2\sqrt{\pi t}}+\frac{3z^2z_{xx}+6z(z_x)^2}{24\sqrt{\pi}t^{3/2}}+\int_{\mathbb{R}^{2}} g_t(\tilde{x})g_t(\tilde{y}) \int_{0}^{h(\tilde{x},\tilde{y})} g_t(z-\tilde{z}) \bigg[\frac{\tilde{x}^2-2t}{4t^2}+\frac{\tilde{x}(\tilde{z}-z)z_x}{2t^2}\\&+\frac{(\tilde{z}-z)^2z_x^2}{4t^2}+\frac{(\tilde{z}-z)z_{xx}}{2t}-\frac{z_x^2}{2t}\bigg] d\tilde{z} d\tilde{x}d\tilde{y}+\text{h.o.t.}
\end{split}
\end{align}
\normalsize
The terms $\frac{z(z_x)^2}{4\sqrt{\pi}t^{3/2}}$ and
\[
\int_{\mathbb{R}^{2}} g_t(\tilde{x})g_t(\tilde{y}) \int_{0}^{h(\tilde{x},\tilde{y})} g_t(z-\tilde{z}) \bigg[\frac{\tilde{x}(\tilde{z}-z)z_x}{2t^2}+\frac{(\tilde{z}-z)^2z_x^2}{4t^2}-\frac{z_x^2}{2t}\bigg] d\tilde{z} d\tilde{x}d\tilde{y}
\]
turn out to be $\mathcal{O}(t^{3/2})$, which is higher than the order needed for the calculations in this paper. We will simplify the two remaining terms in the integrand, starting with the term $\frac{\tilde{x}^2-2t}{4t^2}$. As in the previous calculation substitute in the approximation \cref{eq:intapprox} and integrate:
\small
\begin{align}
\label{eq:dercons}
\begin{split}
&\int_{\mathbb{R}^{2}} g_t(\tilde{x})g_t(\tilde{y}) \int_{0}^{h(\tilde{x},\tilde{y})} g_t(z-\tilde{z}) \frac{\tilde{x}^2-2t}{4t^2} d\tilde{z} d\tilde{x}d\tilde{y}\\
    =&\frac{1}{2\sqrt{\pi t}}\int_{\mathbb{R}^{2}} g_t(\tilde{y}) g_t(\tilde{x}) \Bigg[P[h](\tilde{x},\tilde{y}) -\frac{1}{4t}P[h](\tilde{x},\tilde{y}) z^2\\&+\frac{1}{4t}\big(P[h](\tilde{x},\tilde{y})\big)^2z-\frac{1}{12t}\big(P[h](\tilde{x},\tilde{y})\big)^3\Bigg] \frac{\tilde{x}^2-2t}{4t^2} d\tilde{x}d\tilde{y}+\text{h.o.t.}
\end{split}
\end{align}
\normalsize
Then use \cref{eq:expmom} to further simplify \cref{eq:dercons}:
\small
\begin{align}
\begin{split}
\label{eq:der}
=&\frac{h_{xx}}{2\sqrt{\pi t}}+\frac{\sqrt{t}}{2\sqrt{\pi }}(h_{xxxx}+h_{xxyy})\\
    &-\frac{z^2}{8\sqrt{\pi}t^{3/2}}h_{xx}+\frac{z}{4\sqrt{\pi}t^{3/2}}hh_{xx}+\frac{z}{2\sqrt{\pi t}} \bigg(\frac{3}{2}h_{xx}^2+\frac{1}{2}h_{xx}h_{yy}+h_{xy}^2\bigg) \\&-\frac{h^2}{8\sqrt{\pi }t^{3/2}}h_{xx}-\frac{h}{2\sqrt{\pi t}}\bigg(\frac{3}{2}h_{xx}^2+\frac{1}{2}h_{xx}h_{yy}+h_{xy}^2\bigg)\\&-\frac{\sqrt{t}}{2\sqrt{\pi}}\bigg(\frac{15}{4} h_{xx}^3+\frac{3}{4}h_{xx}h_{yy}^2+\frac{3}{2}h_{xx}^2h_{yy}+6h_{xx}h_{xy}^2+3h_{yy}h_{xy}^2\bigg)+\text{h.o.t}
    \end{split}
\end{align}
\normalsize
We now turn to the $\frac{(\tilde{z}-z)z_{xx}}{2t}$ term, following the same steps:
\small
\begin{align}
\begin{split}
\label{eq:exp2}
&\int_{\mathbb{R}^{2}} g_t(\tilde{x})g_t(\tilde{y}) \int_{0}^{h(\tilde{x},\tilde{y})} g_t(z-\tilde{z}) \frac{(\tilde{z}-z)z_{xx}}{2t} d\tilde{z} d\tilde{x}d\tilde{y}\\
    =&\frac{1}{2\sqrt{\pi t}}\int_{\mathbb{R}^{2}} g_t(\tilde{y}) g_t(\tilde{x}) \Bigg[P[h](\tilde{x},\tilde{y}) -\frac{1}{4t}P[h](\tilde{x},\tilde{y}) z^2\\&+\frac{1}{4t}\big(P[h](\tilde{x},\tilde{y})\big)^2z-\frac{1}{12t}\big(P[h](\tilde{x},\tilde{y})\big)^3\Bigg] \frac{(\tilde{z}-z)z_{xx}}{2t} d\tilde{x}d\tilde{y}+\text{h.o.t.}\\
    =&\frac{z_{xx}}{2\sqrt{\pi t}}\bigg[\frac{h^2}{4t}+\frac{1}{2}h(h_{xx}+h_{yy})+\frac{3t}{4}(h_{xx}^2+h_{yy}^2)+\frac{t}{2}h_{xx}h_{yy}+th_{xy}^2-\frac{zh}{2t}-\frac{z}{2}(h_{xx}+h_{yy})\bigg]\\&+\text{h.o.t.}\\
    \end{split}
\end{align}
\normalsize
Substituting \cref{eq:der} and \cref{eq:exp2} into \cref{eq:conder}, we arrive at the simplification 
\small
\begin{align}
\label{eq:der2}
\begin{split}
\frac{\partial^2}{\partial x^2}& \bigg[G_t*\mathbbm{1}_{\Sigma}\bigg](x,y,z(x,y))|_{(x,y)=(0,0)}\\
=&-\frac{z_{xx}}{2\sqrt{\pi t}}+\frac{z^2z_{xx}}{8\sqrt{\pi}t^{3/2}}+\frac{h_{xx}}{2\sqrt{\pi t}}+\frac{\sqrt{t}}{2\sqrt{\pi}}(h_{xxxx}+h_{xxyy})\\
    &-\frac{z^2}{8\sqrt{\pi}t^{3/2}}h_{xx}+\frac{z}{4\sqrt{\pi}t^{3/2}}hh_{xx}+\frac{z}{2\sqrt{\pi t}} \bigg(\frac{3}{2}h_{xx}^2+\frac{1}{2}h_{xx}h_{yy}+h_{xy}^2\bigg) \\&-\frac{h^2}{8\sqrt{\pi}t^{3/2}}h_{xx}-\frac{h}{2\sqrt{\pi}\sqrt{t}}\bigg(\frac{3}{2}h_{xx}^2+\frac{1}{2}h_{xx}h_{yy}+h_{xy}^2\bigg)\\&-\frac{\sqrt{t}}{2\sqrt{\pi}}\bigg(\frac{15}{4} h_{xx}^3+\frac{3}{4}h_{xx}h_{yy}^2+\frac{3}{2}h_{xx}^2h_{yy}+6h_{xx}h_{xy}^2+3h_{yy}h_{xy}^2\bigg)\\
    &+\frac{z_{xx}}{2\sqrt{\pi t}}\bigg[\frac{h^2}{4t}+\frac{1}{2}h(h_{xx}+h_{yy})+\frac{3t}{4}(h_{xx}^2+h_{yy}^2)+\frac{t}{2}h_{xx}h_{yy}+th_{xy}^2\\&-\frac{zh}{2t}-\frac{z}{2}(h_{xx}+h_{yy})\bigg]+\text{h.o.t.}
\end{split}
\end{align}
\normalsize
Finding 
$\frac{\partial^2}{\partial y^2} \bigg[G_t*\mathbbm{1}_{\Sigma}\bigg](x,y,z(x,y))|_{(x,y)=(0,0)}$ is similar. 
We use both \cref{eq:fconv} and \cref{eq:der2} in several calculations throughout the course of the paper.
\subsection{The $\gamma$s that render the four-step threshold dynamic scheme to be second order accurate}
\label{sec:gamma2}
 We record here the exact values for the coefficients $\gamma$ in the four-stage, second order accurate scheme introduced in \cref{sec:example}.
They are algebraic numbers, but the representations of some of them are quite long and therefore we have approximated them above.
With the exact values given below, we can rigorously state that the \cref{alg:mstage} is second order while maintaining unconditional energy stability. The matrix of values is:
\begin{equation*}
    \gamma=\left(
\begin{array}{cccc}
 1 & 0 & 0 & 0 \\
 -\frac{1}{4} & \frac{5}{4} & 0 &0 \\
 \frac{5}{6} & -\frac{2}{3} & \frac{5}{6} &0 \\
  \gamma_{4,0} & \frac{1}{2} & \gamma_{4,2} &\gamma_{4,3}\\
\end{array}\right).\\
\end{equation*}
\tiny
\begin{align*}
 \gamma_{4,2}=&\bigg(\sqrt[3]{73547857887405865499600064 \sqrt{133495318877644714344377}-23474745371243059566207357648855848671}\\&-\frac{5551049511730043591353151}{\sqrt[3]{73547857887405865499600064 \sqrt{133495318877644714344377}-23474745371243059566207357648855848671}}\\&-456109196575\bigg)/3627134098848
\end{align*}
\begin{align*}
\gamma_{4,3}=&\bigg(-5586815667458\\&\times \sqrt[3]{73547857887405865499600064 \sqrt{133495318877644714344377}-23474745371243059566207357648855848671}\\&+\left(73547857887405865499600064 \sqrt{133495318877644714344377}-23474745371243059566207357648855848671\right)^{2/3}\\&+\frac{31012690382968488487137089701456460158}{\sqrt[3]{73547857887405865499600064 \sqrt{133495318877644714344377}-23474745371243059566207357648855848671}}\\&+\frac{30814150681678355363112149018994128529535197628801}{\left(73547857887405865499600064 \sqrt{133495318877644714344377}-23474745371243059566207357648855848671\right)^{2/3}}\\&+7974522440634228925392639\bigg)/18386964471851466374900016
\end{align*}
\normalsize
\small
\[
\gamma_{4,0}=1-\frac{1}{2}-\gamma_{4,2}-\gamma_{4,3}
\]
\normalsize
It can be checked that these $\gamma$s satisfy the inequalities in the hypothesis of \cref{claim:ms} for stability, and the consistency equations in \cref{claim:cons} for second order exactly.

\section*{Acknowledgments}
Funding: This work was supported by National Science Foundation [grant numbers DMS-1719727, DMREF-1436154, DMREF-1729166].

\bibliographystyle{siamplain}
\bibliography{references}
\end{document}